\documentclass[12pt]{amsart}
\usepackage{amssymb, amsmath}
\usepackage{latexsym}
\usepackage{citesort}
\usepackage{epsfig}

\newdimen\AAdi%
\newbox\AAbo%
\def\AAk#1#2{\s_etbox\AAbo=\hbox{#2}\AAdi=\wd\AAbo\kern#1\AAdi{}}%
\def\AAr#1#2#3{\s_etbox\AAbo=\hbox{#2}\AAdi=\ht\AAbo\raise#1\AAdi\hbox{#3}}%
\font\tenmsb=msbm10 at 12pt
\font\sevenmsb=msbm7 at 8pt
\font\fivemsb=msbm5 at 6pt
\newfam\msbfam
\textfont\msbfam=\tenmsb
\scriptfont\msbfam=\sevenmsb
\scriptscriptfont\msbfam=\fivemsb

\newtheorem{thm}{Theorem}
\newtheorem{lem}{Lemma}
\newtheorem{cor}{Corollary}
\newtheorem{rem}{Remark}
\newtheorem{pro}{Proposition}

\newcommand{\ba}{\begin{array}}
\newcommand{\ea}{\end{array}}

\def\beq{\begin{equation}}
\def\eeq{\end{equation}}

\def\real     #1{{\mathbb R^{#1}}}

\def\d        #1#2{{#1}_{#2}}
\def\dd       #1#2#3{{#1}_{#2#3}}

\def\dddd     #1#2#3#4#5{{#1}_{#2#3#4#5}}

\def\uu       #1#2#3{{#1}^{#2#3}}

\def\uddd     #1#2#3#4#5{#1^#2_{\phantom{#2}#3#4#5}}

\def\uudd     #1#2#3#4#5{#1^{#2#3}_{\phantom{#2#3}#4#5}}
\def\udud     #1#2#3#4#5{#1^{#2\phantom{#3}#4}_
              {\phantom{#2}#3\phantom{#4}#5}}

\def\dduu     #1#2#3#4#5{#1_{#2#3}^{\phantom{#2#3}#4#5}}
\def\dudu     #1#2#3#4#5{#1_{#2\phantom{#3}#4}^
              {\phantom{#2}#3\phantom{#4}#5}}
\def\duud     #1#2#3#4#5{#1_{#2\phantom{#3#4}#5}^{\phantom{#2}#3#4}}

 \def \vs{\vspace*{0.2cm}}

\def\ds{\displaystyle}
\def\p{\partial}

\def\b{\beta}

\def\l{\lambda}

\def\D{\Delta}
\def\d{\delta}
\def\s{\sigma}

\def\a{\alpha}

\def\cal{\mathcal}

\def\M1{{{\cal M}}_1}

\def\n{\nabla}

\def\R{\mathbb R}

\def\p{\partial}
\def\dt{\frac {d}{dt}}

\begin{document}
\title
[Mean curvature flow]
{Mean curvature flow with flat normal bundles}

\author[K. Smoczyk]{Knut Smoczyk}
\address[Knut Smoczyk]{Max Planck Institute for Mathematics in the Sciences\\
                       Inselstr.\ 22-26\\
                       04103 Leipzig, Germany}
\email{smoczyk@mis.mpg.de}
\thanks{The research of the first author was supported by a Heisenberg
fellowship of the DFG}

\author[G. Wang]{Guofang Wang}
\address[Guofang Wang]{Institute of Mathematics\\
                        Fudan University\\
                        Shanghai, China
                        and
                        Max Planck Institute for Mathematics in the Sciences\\
                       Inselstr.\ 22-26\\
                       04103 Leipzig, Germany}
\email{gwang@mis.mpg.de}
%\thanks{}

\author[Y. L. Xin]{Y.~L.~ Xin}
\address[Y.~L.~ Xin]{Institute of Mathematics\\
                        Fudan University\\
                        Shanghai, China}
\email{ylxin@fudan.edu.cn}
\thanks{The research of the third author was partially supported by project
\# 973 of MSTC and SFECC}

\begin{abstract}
We show that flatness of the normal bundle is preserved under the
mean curvature flow in $\real{n}$ and use this to
generalize a classical result for hypersurfaces
due to Ecker \& Huisken \cite{ecker huisken} in the case of submanifolds
with arbitrary codimension.
\end{abstract}
\renewcommand{\subjclassname}{%
  \textup{2000} Mathematics Subject Classification}
\date{Oct. 31, 2004}%\today}
\maketitle

%%%%%%%%%%%%%%%%%%%%%%%%%%%%%%%%%%%%%%%%%%%%%%%%%%%%%%%%%%%%%%%%%%%%%%%%%%%%%%
\section{Introduction}
Let us consider immersions
$$F:M^m\to\real{n}$$
of an $m$-dimensional submanifold in $\real{n}$ of codimension $k$.
Throughout this paper we shall assume that there is a one-parameter
family $F_t=F(\cdot,t)$ of immersions with corresponding images
$M_t=F_t(M)$ such that {\sl mean curvature flow}
\begin{eqnarray}
\dt F(x,t)&=&H(x,t)\,,\quad x\,\in\,M\label{mcf}\\
F(x,0)&=&F_0(x)\nonumber
\end{eqnarray}
is satisfied for some initial data $F_0$.
Here, $H(x,t)$ is the mean curvature vector of $M_t$ at the point
$x\in M$, i.e. $H$ is the trace of the second fundamental form
\begin{equation}\nonumber
A=\nabla dF.
\end{equation}
The mean curvature flow has been studied intensively by many authors.
Most of the results have been obtained for hypersurfaces. A classical
result is due to Ecker and Huisken \cite{ecker huisken},
where they study hypersurfaces
in $\real{n+1}$ that can be represented as entire graphs over a flat
plane. Their result says that any polynomial growth rate for the
height and the gradient of the initial surface $M_0$ is preserved during the
evolution and that in case of Lipschitz initial data with linear growth,
(\ref{mcf}) has a smooth solution for all times $t>0$.
The growth condition was
$$v=\frac{1}{\langle\nu,\omega\rangle}\le \text{const},$$
where $\omega\in\real{n+1}$, $|\omega|=1$ and $\nu$ is some choice
of unit normal vector such that $\langle\nu,\omega\rangle>0$.
In addition, they proved
that these hypersurfaces approach a {\sl selfsimilar expanding} solution
of (\ref{mcf}) as $t\to\infty$ provided the initial graph was ``straight''
at infinity.

In higher codimension, the general picture of mean curvature
flow is still incomplete even though some work has been carried
out by Chen, Li \cite{chen li}, Wang \cite{wang}, \cite{wang2} and one of
the authors \cite{smoczyk wang}. E.g. based on the interior
estimates for hypersurfaces obtained by Ecker and Huisken in
\cite{ecker huisken2}, Wang \cite{wang2}
proved that any initial compact submanifold
that satisfies a local K-Lipschitz condition admits a smooth solution
on some time interval $(0,T), T>0$.

In the case of graphical mean curvature flow in higher codimension
there are some longtime existence and convergence results to flat
spaces. In \cite{wang} Wang defined a similar expression as
the above mentioned quantity $v$ which was essential to get the
a-priori estimates needed in the longtime existence and convergence
results. However, in that paper the author had to assume a smallness
condition on $v$ and the theorem did not apply to arbitrary graphs.
In \cite{smoczyk wang} it was shown that
Lagrangian graphs with convex potentials admit smooth solutions
for all times and that the solutions exponentially converge to
flat Lagrangian planes.

So, for some time it was unclear how to generalize the results in
\cite{ecker huisken} in the ``best'' way to the case of arbitrary
codimension. Motivated by a recent paper of the third author \cite{xin},
we believe that the flatness of the normal bundle is the key
ingredient to get the convergence results for arbitrary graphs
in higher codimension. Note, that trivially the normal bundle of
any hypersurface is flat. So our result is a natural extension
of the results by Ecker and Huisken.

There are many examples of submanifolds with flat normal bundles.
Hypersurfaces are trivial examples.
Other trivial examples are curves in $\real{n}$.
In addition, any submanifold of codimension
$2$ in $\real{n}$ which also is a hypersurface of the
standard sphere $S^{n-1}$
% or any other hypersurface of constant mean curvature
must have a flat normal bundle. For more examples, see
\cite{hpt}, \cite{terng1}, \cite{terng2} and \cite{terng3}, where a theory
of isoparametric submanifolds
was established in the framework of flat normal bundles.

The organization of the paper is as follows: In section \ref{class}
we recall the monotonicity formula and the noncompact maximum principle
by Ecker and Huisken and introduce the class of solutions for which
our results will apply. Section \ref{geometric} introduces our notation
and recalls the most important structure equations in the geometry
of real submanifolds. Some basic evolution equations for the mean
curvature flow are derived resp. recalled in section \ref{evolution}.
In section \ref{growth} we prove that polynomial growth rates
are preserved in arbitrary codimension and that these growth
estimates can be applied even to non-graphical submanifolds, like
cylinders. The core of our article is the proof of Theorem
\ref{flatness}, i.e. that flatness of the normal bundle is preserved.
This will be done in section \ref{flatness 2}. In the remainder -
based on this fundamental observation - we can proceed basically as in
\cite{ecker huisken}, to carry over the results by Ecker and Huisken
to the case of arbitrary codimension.

We are indebted to J\"urgen Jost for his constant support and the MPI in
Leipzig for hospitality. The authors also wish to thank Klaus Ecker for
fruitful discussions.

%%%%%%%%%%%%%%%%%%%%%%%%%%%%%%%%%%%%%%%%%%%%%%%%%%%%%%%%%%%%%%%%%%%%%%%%%%%%%%
\section{The class of solutions}\label{class}
Throughout this article we shall assume that the solutions
of (\ref{mcf}) that we consider will belong to the class of solutions
for which we can apply Huisken's {\sl monotonicity formula}
\begin{pro}(Huisken)
For a function $f(x,t)$ on $M$ we have
\begin{equation}\label{mon}
\dt\int_M f\rho d\mu_t=\int_M\left(\dt\,f-\Delta f\right)\rho d\mu_t
-\int_M f\rho\left|H+\frac{F^\perp}{2(t_0-t)}\right|^2d\mu_t,
\end{equation}
where
$$\rho(y,t)=\frac{1}{(4\pi(t_0-t))^\frac{n}{2}}\,e^{-\frac{|y|^2}{4(t_0-t)}}$$
is the {\sl backward heat kernel} on $\real{n}$ at the origin, and $t_0>t$.
\end{pro}
In particular, we will assume that integration
by parts is permitted and all integrals are finite
for the submanifolds and functions we will consider
in the sequel. E.g., this is the case for those smooth solutions of
(\ref{mcf}) for
which the curvature and its covariant derivatives
have at most polynomial growth
at infinity since then the faster exponential decay rate of the heat kernel
$\rho$ yields finite integrals.

As in \cite{ecker huisken} we will repeatedly use the following
{\sl maximum principle} which is based on the monotonicity formula.

\begin{cor}(Ecker and Huisken \cite{ecker huisken})\label{max princ}
Suppose the function $f=f(x,t)$ satisfies the inequality
\begin{equation}\nonumber
\left(\dt-\Delta\right)f\le\langle {\bf a},\nabla f\rangle
\end{equation}
for some $\bf a$
which is uniformly bounded on $M\times[0,t_1]$ for some $t_1>0$, then
$$\sup_{M_t}f\le\sup_{M_0}f$$
for all $t\in[0,t_1]$.
\end{cor}

%%%%%%%%%%%%%%%%%%%%%%%%%%%%%%%%%%%%%%%%%%%%%%%%%%%%%%%%%%%%%%%%%%%%%%%%%%%%%%
\section{Geometric quantities of immersions}\label{geometric}
Let $F:M^m\to \R^{n}$ be an immersion and let $k$ be the codimension of
$M$, i.e. $n=k+m$.
We let $(x^i)_{i=1,\dots,m}$ denote local coordinates on $M$ and we will
always use {\sl cartesian} coordinates $(y^\alpha)_{\alpha=1,\dots,n}$
on $\real{n}$. Doubled greek and latin indices are summed form $1$ to $n$
resp. from $1$ to $m$.
In local coordinates the {\sl differential} $dF$ of $F$ is given by
$$dF=F^\alpha_i\frac{\partial}{\partial y^\alpha}\otimes dx^i,$$
where $F^\alpha=y^\alpha(F)$ and
$F^\alpha_i=\frac{\partial F^\alpha}{\partial x^i}$.
The coefficients of the induced metric $g_{ij}\,dx^i\otimes dx^j$ are
\[g_{ij}=\langle F_i,F_j\rangle=g_{\alpha\beta}F^\alpha_iF^\beta_j,\]
where $g_{\alpha\beta}=\delta_{\alpha\beta}$ is the euclidean metric
in cartesian coordinates. As usual, the Christoffel symbols are
\[\Gamma^k_{ij}=\frac 12 g^{kl}(\frac{\partial g_{lj}}{\partial x^i}+
\frac{\partial g_{li}}{\partial x^j}-
\frac{\partial g_{ij}}{\partial x^l}).\]

The second fundmental form is defined by
\[A=\n dF:=A^\a_{ij}\,\frac{\partial }{\partial y^\a}\otimes
dx^i\otimes dx^j.\]
Here and in the following all canonically induced full connections on
bundles over $M$ will be denoted by $\nabla$. Later, we will occasionally
also use the connection on the normal bundle which will then be denoted
by $\nabla^\perp$.
It is easy to check that in cartesian coordinates on $\real{n}$ we have
\[A^\a_{ij}=F^\a_{ij}-\Gamma^k_{ij}F^\a_k,\]
where $F^\alpha_{ij}=\frac{\partial^2 F^\alpha}{\partial x^i\partial x^j}$.
By definition, $A$ is a section in $F^{-1}T\real{n}\otimes T^*M\otimes T^*M$
and it can be easily checked that $A$ is normal, i.e. that
$$A\in\Gamma\left(NM\otimes T^*M\otimes T^*M\right),$$
where $NM$ denotes the normal bundle of $M$ w.r.t. the immersion $F$.
This means that
\begin{equation}\label{normal}
g_{\alpha\beta}A^\alpha_{ij}F^\beta_k=0,\quad\forall\,i,j,k.
\end{equation}
In particular, the {\sl mean curvature vector field}
$H=H^\alpha\frac{\partial}{\partial y^\alpha}$ with
$H^\alpha=g^{ij}A^\a_{ij}$ satisfies
\[g_{\a\b}H^\a F^\b_j=0,\quad\forall\,j.\]

The curvature of the normal bundle is defined locally by
$\uudd R\alpha\beta ij
\frac{\p}{\p y^\a}\otimes \frac{\p}{\p y^\b}\otimes dx^i\otimes dx^j$,
where by {\sl Ricci's equation}
\begin{equation}\label{ricci}
\uudd R\alpha\beta ij= A^\a_{is}A^{\b s}_{j}-A^\a_{js}A^{\b s}_{i}.
\end{equation}
The normal  bundle is flat if and only if
$\uudd R\alpha\beta ij$ vanishes for any $\a,$ $\b$, $i$ and $j$.
In addition we have the {\sl Gau\ss\ equations} for the induced curvature
tensor on $M$
\begin{equation}\label{gauss}
\dddd Rijkl=
\dd g\alpha\beta(A^\alpha_{ik}A^\beta_{jl}-A^\alpha_{il}A^\beta_{jk}).
\end{equation}
Let us also recall the {\sl Codazzi equations}
\begin{equation}\label{1} \n_i A^\a_{jk} =\n _jA^\a_{ik}+F^\a_l\uddd Rlkji,
\end{equation}

\begin{equation}\label{2} \n^k A^\a_{jk} =\n_jH^\a+F^\a_lR^l_{j}.
\end{equation}

The following rule for interchanging derivatives
\begin{equation}\label{3}
\n_i\n_jA^\a_{lk}=\n_j\n_i A^\a_{lk}-\uddd RmlijA^\a_{mk}-\uddd RmkijA^\a_{lm}
\end{equation}
and the second Bianchi identity together with the Codazzi equations imply
the {\sl Simons' identity}

\begin{equation}\label{4}\ba{rcl}
\ds\vs \D A^\a_{lk} &=& \ds \n_l\n_k H^\a +R^m_lA^\a_{mk}+R^m_kA^\a_{ml}
-2A^\a_{jm}\dudu Rljkm\\
&&\ds+F^\a_m(\n_lR^m_k+\n_kR^m_l-\n^mR_{lk}),\ea
\end{equation}
where $\dd Rij$ denotes the Ricci curvature of $M$.

%%%%%%%%%%%%%%%%%%%%%%%%%%%%%%%%%%%%%%%%%%%%%%%%%%%%%%%%%%%%%%%%%%%%%%%%%%%%%%
\section{Evolution equations}\label{evolution}
>From the main evolution equation
\[\dt F=H\]
we obtain
\[\begin{array}{rcl}
\ds\vs\dt F^\a_i &=& \ds \nabla_iH^\a.
\end{array}\]
Let us define the symmetric tensors
$$\dd aij:=\dd g\alpha\beta H^\alpha A^\beta_{ij}\,,\quad
\dd bij:=\dd g\alpha\beta A^\alpha_{ik}A^{\beta k}_j,$$
so that by Gau\ss' equation the Ricci tensor satisfies
\begin{equation}\label{Rij}
\dd Rij=\dd aij-\dd bij\,.
\end{equation}
>From $\dd g\alpha\beta F^\alpha_i H^\beta=0$ we derive
the evolution equation for the induced metric $g_{ij}$
\begin{equation}\label{metric}
\dt g_{ij}=2g_{\a\b}\nabla_iH^\a F^\b_j=-2g_{\a\b}H^\a A^\b_{ij}=-2\dd aij.
\end{equation}
Consequently, the volume form $d\mu$ on $M$ satisfies
\begin{equation}\label{volume}
\dt\,d\mu=-|H|^2d\mu.
\end{equation}

We need to compute the evolution equation for $A^\a_{ij}$. In a first
step we get
\[\begin{array}{rcl}
\ds\vs \dt A^a_{ij}&=& \ds \dt F^\a_{ij}- F^\a_k \dt \Gamma ^k_{ij}-
\Gamma^k_{ij}
\dt F^\a_k.\\
&=&\ds \n_i\n_j H^\a-\dt \Gamma^k_{ij} F^\a_k.\end{array}\]
The evolution equation for $\Gamma^\a_{ij}$ is
\[\dt \Gamma^k_{ij}=\frac 12 g^{kl} (\n_i \dt g_{lj}+\n_j \dt g_{li}-
\n_{l} \dt g_{ij}).\]

The last two equations, the Simons' identity and (\ref{Rij}), (\ref{metric})
imply
\begin{eqnarray}
\dt A^\a_{ij} &=& \ds \D A^\a_{ij}-R^m_i A^\a_{mj}-R^m_j A^\a_{mi}+
2A^\a_{mn} \udud Rminj\nonumber\\
&&-F^\a_m(\n_i b^m_j+\n_j b^m_i-\n^m b_{ij})\label{A}
\end{eqnarray}

In addition
\begin{eqnarray}
\dt H^\alpha&=&\uu gij\dt A^\alpha_{ij}+2\uu aij A^\alpha_{ij}\nonumber\\
&=&\uu gij\left(\nabla_i\nabla_jH^\alpha-\dt\Gamma_{ij}^k F^\alpha_k\right)
+2\uu aij A^\alpha_{ij}\nonumber\\
&=&\Delta H^\alpha-\uu gij\dt\Gamma_{ij}^k F^\alpha_k+2\uu aijA^\alpha_{ij}
\nonumber
\end{eqnarray}

so that
\begin{eqnarray}
\dt |H|^2&=&2H_\alpha\dt H^\alpha\nonumber\\
&=&2H_\alpha\Delta H^\alpha+4|a_{ij}|^2,\nonumber
\end{eqnarray}
for $F^\a_kH_\a=0$.
Thus
\begin{equation}\nonumber
\dt |H|^2= \D |H|^2-2|\n_i H^\a|^2+4|\dd aij|^2.
\end{equation}
Since $H_\alpha F^\alpha_i=0$ we conclude that
$$\nabla_iH_\beta F^\beta_jF^\alpha_l\uu gjl=-a_i^lF^\alpha_l$$
so that
$$|\nabla_iH^\alpha|^2=|\nabla_iH^\alpha+a_i^lF_l^\alpha|^2+|\dd aij|^2.$$
Hence
\begin{equation}\label{mean}
\dt |H|^2= \D |H|^2-2|\n_i H^\a+a_i^lF^\alpha_l|^2+2|\dd aij|^2.
\end{equation}

Let $\n^\perp$ denote the normal connection
induced from the immersion $F$.
Then $|\n^\perp H|^2=|\n_i H^\a+a_i^lF^\alpha_l|^2$.
\begin{rem}
For a hypersurface with inward unit normal vector $\nu$, scalar mean
curvature $H$ and second fundamental tensor $\dd hij$ we have
$A^\alpha_{ij}=\dd hij\nu^\alpha, H^\alpha=H\nu^\alpha$ and
$\dd aij=H\dd hij, \nabla_iH^\alpha=\nabla_iH\nu^\alpha-a_i^lF_l^\alpha$
so that $-2|\nabla_iH^\alpha|^2=-2|\nabla H|^2-2|\dd aij|^2$ and
$|\dd aij|^2=H^2|A|^2$.
\end{rem}

Now we compute the evolution equation for $|A|^2$. From (\ref{metric}),
(\ref{A}) and the normality of $A$ we deduce
\begin{eqnarray}
\dt |A|^2 &=&4\uu aij\dd bij+2A_\alpha^{ij}\dt A^\alpha_{ij}\nonumber\\
&=&4\uu aij\dd bij+2A_\alpha^{ij}\left(\Delta A^\alpha_{ij}
-R^m_i A^\a_{mj}-R^m_j A^\a_{mi}+
2A^\a_{mn} \udud Rminj\right)\nonumber\\
&=&\Delta|A|^2-2|\nabla_lA^\alpha_{ij}|^2+4|\dd bij|^2+4A_\alpha^{ij}
A^\alpha_{mn}\udud Rminj.
\end{eqnarray}
Since
\begin{eqnarray}
4|\dd bij|^2+4A_\alpha^{ij}A^\alpha_{mn}\udud Rminj
&=&2|\uudd R\alpha\beta ij|^2+4|A^{\alpha m}_nA^{\beta n}_m|^2\nonumber
\end{eqnarray}
and the tangential part of $\nabla_lA^\alpha_{ij}$ is given by
$-A^\beta_{ij}A^k_{\beta l}F^\alpha_k$ we conclude
$$|\nabla_lA^\alpha_{ij}|^2=|\nabla_lA^\alpha_{ij}
+A^\beta_{ij}A^k_{\beta l}F^\alpha_k|^2+|A^{\alpha m}_nA^{\beta n}_m|^2$$
and finally
\begin{equation}\label{A2}
\dt|A|^2=\Delta|A|^2
-2|\nabla_lA^\alpha_{ij}+A^\beta_{ij}A^k_{\beta l}F^\alpha_k|^2
+2|A^{\alpha m}_nA^{\beta n}_m|^2+2|\uudd R\alpha\beta ij|^2.
\end{equation}

Here $|\nabla_lA^\alpha_{ij}+A^\beta_{ij}A^k_{\beta l}F^\alpha_k|^2=|\n^\perp A|^2$.
%%%%%%%%%%%%%%%%%%%%%%%%%%%%%%%%%%%%%%%%%%%%%%%%%%%%%%%%%%%%%%%%%%%%%%%%%%%%%%
\section{Growth estimates}\label{growth}
\begin{figure}[ht]
\begin{center}
\epsfig{file=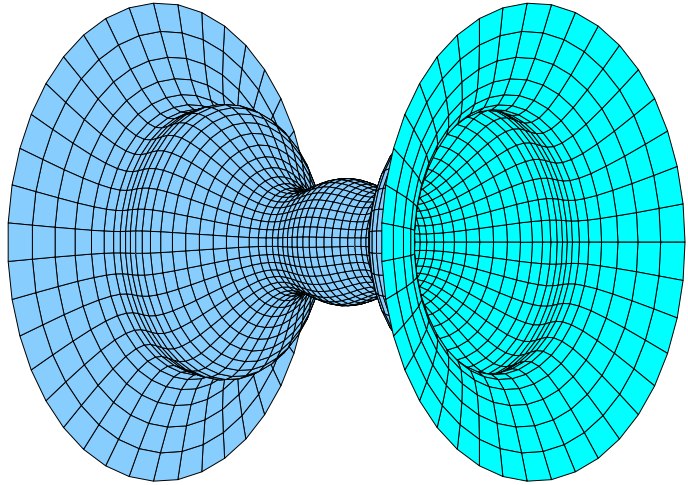,width=.8\hsize}
\caption{\sl Growth rates of cylindrical surfaces are preserved.
Here, the surface is given by
$
(x, \phi)\mapsto (u_1,u_2,u_3)=
\bigl(x,(1+|x|)\cos x\cos\phi,(1+|x|)\cos x\sin\phi\bigr).
$
Since $u^2:=u_2^2+u_3^2=(1+|x|)^2\cos^2x\le 2(1+x^2)$ the surface
grows linearly over the axis of rotation. We have
$\ell=1$, $m=2$ and $n=3$.}
\label{figure 1}
\end{center}
\end{figure}

In this section we will prove that polynomial growth estimates for
submanifolds in $\real{n}$ are preserved under the mean curvature flow.
Note, that in this section we do not require that $M_0$ or $M_t$ can
be written as graphs over some flat $m$-plane. We will use a
flat $\ell$-plane merely as a reference submanifold to measure distances.
As a consequence, we will also obtain growth estimates for other objects
as graphs, e.g. for cylindrical objects as depicted in figure \ref{figure 1}.

Choose an orthonormal basis $e_1,\dots, e_{n}$ of $\real{n}$
and let $1\le\ell\le n$ be an integer (for later purposes we have in mind
$\ell=m$ but here in general $\ell$ and $m$ can be different).
For an immersion $F:M^m\to\real{n}$ we define $n$ {\sl coordinate functions}
$$u_i:=\langle F,e_i\rangle\,,\quad i=1,\dots,n$$
In addition we define
$$x:=\sqrt{\sum_{i=1}^\ell u_i^2}\,,\quad u:=\sqrt{\sum_{i=\ell+1}^n u_i^2}$$
such that $|F|^2=x^2+u^2$. Since $\left(\dt-\Delta\right)u_i=0$ we conclude
\begin{eqnarray}
\dt |F|^2&=&\Delta |F|^2-2m,\label{eq 1}\\
\dt u^2&=&\Delta u^2-2\sum_{i=\ell+1}^n|\nabla u_i|^2,\label{eq 2}\\
\dt x^2&=&\Delta x^2-2m+2\sum_{i=\ell+1}^n|\nabla u_i|^2,\label{eq 3}
\end{eqnarray}
since
$$\sum_{i=1}^n|\nabla u_i|^2=m.$$
For a constant $c>0$ to be determined later, we define the function
$$\eta:=1+ct+x^2.$$
Let $\varphi:\real{+}\to\real{+}$ be a smooth function with
$\varphi'\le0$, $\varphi''\ge 0$.
We want to compute the evolution equation of $u^2\varphi(\eta)$.
\begin{eqnarray}
\dt\left(u^2\varphi\right)
&=&\varphi\left(\Delta u^2-2\sum_{i=\ell+1}^n|\nabla u_i|^2\right)
+u^2\varphi'\left(\Delta\eta+2\sum_{i=\ell+1}^n|\nabla u_i|^2
+c-2m\right)\nonumber\\
&=&\Delta\left(u^2\varphi\right)-2\varphi'\langle\nabla\eta,\nabla u^2\rangle
-\varphi''u^2|\nabla\eta|^2
-2(\varphi-u^2\varphi')\sum_{i=\ell+1}^n|\nabla u_i|^2
\nonumber\\
&&+(c-2m)u^2\varphi'.\nonumber
\end{eqnarray}
If at some point on $M_t$ we have $u^2=0$, then  at such a point
$$\dt\left(u^2\varphi\right)\le\Delta\left(u^2\varphi\right).$$
We want to prove that this inequality holds at all points on $M_t$.
At those points, where $\varphi'=0$ we are done as well. So w.l.o.g. we
can assume that $u\neq0$ and $\varphi'\neq0$.
Next observe that
$$|\nabla u^2|^2\le 4u^2\sum_{i=\ell+1}^n|\nabla u_i|^2$$
implies
$$|\nabla u|^2\le \sum_{i=\ell+1}^n|\nabla u_i|^2$$
at all points, where $u\neq 0$.
We use Schwarz' inequality to estimate
$$-2\varphi'\langle\nabla\eta,\nabla u^2\rangle
\le -2\varphi'(\varepsilon
u^2|\nabla\eta|^2+\frac{1}{\varepsilon}|\nabla u|^2)
\le -2\varepsilon\varphi'u^2|\nabla\eta|^2
-\frac{2}{\varepsilon}\,\varphi'\sum_{i=\ell+1}^n|\nabla u_i|^2,$$
where $\varepsilon$ is some positive constant. We choose $\varepsilon=
-\frac{\varphi'}{\varphi}$ and get
\begin{eqnarray}
\dt\left(u^2\varphi\right)
&\le&\Delta\left(u^2\varphi\right)
+\left(\frac{2(\varphi')^2}{\varphi}-\varphi''\right)
u^2|\nabla\eta|^2+(c-2m)u^2\varphi'.\nonumber
\end{eqnarray}
As in \cite{ecker huisken} we have
$$\nabla_l\eta=2\left\langle F_l,F-\sum_{i=\ell+1}^nu_ie_i\right\rangle,$$
and therefore
$$|\nabla\eta|^2\le 4|F-\sum_{i=\ell+1}^nu_ie_i|^2=4x^2\le 4\eta,$$
so that always
\begin{eqnarray}
\dt\left(u^2\varphi\right)
&\le&\Delta\left(u^2\varphi\right)
+u^2\left\{
4\eta\left(\frac{2(\varphi')^2}{\varphi}-\varphi''\right)
+(c-2m)\varphi'\right\}.\label{u2p}
\end{eqnarray}
\begin{pro}\label{pro0}
If for some $c_0<\infty, p\ge 0$, the inequality
$$u^2\le c_0\left(1+|F|^2-u^2\right)^p$$
is satisfied on $M_0$, then for all $t>0$,
$$u^2\le c_0\left(1+|F|^2-u^2
+\bigl(2m+4(p-1)\bigr)t\right)^p.$$
\end{pro}
\begin{proof}
We choose $\varphi(\eta)=\eta^{-p}$ and $c=2m+4(p-1)$.
Then $\varphi'=-p\eta^{-p-1}$,
$\varphi''=p(p+1)\eta^{-p-2}$. Inserting this into (\ref{u2p})
gives
\begin{eqnarray}
\dt\left(u^2\varphi\right)
&\le&\Delta\left(u^2\varphi\right)
+\frac{u^2}{\eta}\left(8p^2-4p(p+1)-p(c-2m)\right)\nonumber\\
&=&\Delta\left(u^2\varphi\right).\nonumber
\end{eqnarray}
and the result follows from Corollary \ref{max princ}.
\end{proof}

%%%%%%%%%%%%%%%%%%%%%%%%%%%%%%%%%%%%%%%%%%%%%%%%%%%%%%%%%%%%%%%%%%%%%%%%%%%%%%
\section{Preserving flatness of the normal bundle}\label{flatness 2}
In this section we will prove that flatness of the normal bundle is preserved.
We do not require that $M$ is a graph nor do we assume compactness or
completeness. The theorem can be applied to any smooth solution of the
mean curvature flow for which $M_0$ has a flat normal bundle.
\begin{thm}\label{flatness}
Let $F:M\times[0,T)\to\real{n}$ be a smooth solution of the mean curvature
flow and assume that the normal bundle of $M_0$ is flat. If $|A|^2$ is bounded
on each $M_t$ then the normal bundle of $M_t$ is flat as well.
\end{thm}

Note that we do not require that $|A|^2$ is uniformly bounded in $t$.

\begin{proof}
For the proof of this theorem we have to compute the evolution equation
of the normal curvature tensor $\uudd R\alpha\beta ij$. We will show
that the squared normal curvature tensor $R^\perp$ satisfies an evolution
equation of the form

\begin{equation}\label{normal 2}
\dt|R^\perp|^2=\Delta|R^\perp|^2-2|\nabla R^\perp|^2+A*A*R^\perp*R^\perp,
\end{equation}

where the last term is a contraction of a term quadratic in $A$ and
one which is quadratic in $R^\perp$.
Then, by assumption on $|A|^2$, on a compact time interval $[0,t_1]$
we can choose a constant $c$ (depending on $t_1$) such that

$$\dt |R^\perp|^2\le\Delta|R^\perp|^2-2|\nabla R^\perp|^2+c|R^\perp|^2$$

and the function $f:=e^{-ct}|R^\perp|^2$ satisfies
$$\dt f\le \Delta f$$
on $[0,t_1]$. The result then follows from Corollary \ref{max princ}.

It remains to derive the evolution equation for $|R^\perp|^2$.
It turns out
that the computation is rather complicated, since a number of
symmetries of the curvature tensor and the second fundamental form
have to be used.

The first observation is
\[|R^\perp|^2=|\uudd R\alpha\beta ij|^2=2|b_{ij}|^2-2c^{\a\b}_{ij}c^{ij}_{\b\a},\]
where \[c^{\a\b}_{ij}
=A^\a_{ik}A^{\b k}_j.\]

\begin{eqnarray}
\ds \dt b_{ij}
& =&\ds\vs 2 a^{kl}A_{\a ik}A^\a_{lj}+ \dt A^\a_{kj} A^k_{\a i}
+ \dt A^\a_{ki} A^k_{\a j}\nonumber\\
&=&\ds\vs  2 a^{kl}A_{\a ik}A^\a_{lj}\nonumber\\
&&\ds\vs + A^k_{\a j}
(\D A^\a_{ik}-R^m_iA^\a _{mk}
-R^m_k A^\a _{mi}+2A^\a_{mn}{{{R^m}_i}^n}_k)\nonumber\\
&&\ds\vs + A^k_{\a i}
(\D A^\a_{jk}-R^m_jA^\a _{mk}
-R^m_k A^\a _{mj}+2A^\a_{mn}{{{R^m}_j}^n}_k)\nonumber
\end{eqnarray}
and then
\begin{eqnarray}
\ds \dt b_{ij}
&=&\ds\vs  2 a^{kl}A_{\a ik}A^\a_{lj}+ \D b_{ij}\nonumber\\
&&\ds\vs -\n^l A^k_{\a j} \n_l A^\a_{ik}-
R^m_ib_{jm}-R^m_k A^k_{\a j}A^\a_{mi}
+2A^k_{\a j}A^\a_{mn}{{{R^m}_i}^n}_k\nonumber\\
&&\ds\vs -\n^l A^k_{\a i} \n_l A^\a_{jk}-
R^m_jb_{im}-R^m_k A^k_{\a i}A^\a_{mj}
+2A^k_{\a i}A^\a_{mn}{{{R^m}_j}^n}_k.\nonumber
\end{eqnarray}

This implies
\begin{eqnarray}
\ds\dt |b_{ij}|^2 &=&\ds \vs 4a^{ik} b_{il} b^l_k+ 2b^{ij}\dt b_{ij}
\nonumber\\
&=&\ds\vs 4a^{ik}b_{il}b^l_k+ 4a^{kl}A^\a_{ik}A_{\a lj}b^{ij}
+2b^{ij}(\D b_{ij}-2\n^l A^k_{\a j} \n_l A^\a_{ik})\nonumber\\
&&\ds\vs +4 b^{ij}(-R^m_ib_{jm}-R^m_k A^k_{\a j} A^\a_{mi}
+2 A^k_{\a j}A^\a_{mn} {{{R^m}_i}^n}_k)\nonumber
\end{eqnarray}
and from $\dd Rij=\dd aij-\dd bij$
\begin{eqnarray}
\ds\dt |b_{ij}|^2
&=&\ds\vs 4b^{ij}b^m_ib_{jm}+4b^{ij}b^m_k A^k_{\a j} A^\a_{mi}\nonumber\\
&&\ds\vs +\D|b_{ij}|^2-2|\n_l b_{ij}|^2-4b^{ji}\n^lA^k_{\a j}\n_l A^\a_{ik}
\nonumber\\
&&\ds\vs +8b^{ij}A^k_{\a j}A^\a_{mn} {{{R^m}_i}^n}_k\nonumber\\
&=&\ds\vs 4b^{ij}b^m_ib_{jm}+4b^{ij}b^m_k A^k_{\a j} A^\a_{mi}\nonumber\\
&&\ds\vs +\D|b_{ij}|^2-2|\n_l b_{ij}|^2-4b^{ji}\n^lA^k_{\a j}\n_l A^\a_{ik}
\nonumber\\
&&\ds +8b^{ij}A^k_{\a j}A^{\a mn}(A^\b_{mn}A_{\b ik}-A^\b_{mk}A_{\b in})
\nonumber\\
&=&\ds\vs \D|b_{ij}|^2-2|\n_l b_{ij}|^2-4b^{ji}\n^lA^k_{\a j}\n_l A^\a_{ik}
\nonumber\\
&&\ds\vs +4b^{ij}b^m_ib_{jm}+4b^{ij}b^m_k A^k_{\a j} A^\a_{mi}
+8b^{ij}g_{ln}c^{\a\b}_{ji}c^{nl}_{\a\b}-
8b^j_lc^{\a\b}_{jk}c^{lk}_{\b\a}\nonumber
\end{eqnarray}
so that
\begin{eqnarray}
\dt|\dd bij|^2
&=&\ds\vs \D|b_{ij}|^2-2|\n_l b_{ij}|^2-4b^{ji}\n^lA^k_{\a j}\n_l A^\a_{ik}
\ds\vs +4\Gamma_1+4\Gamma_2+8\Gamma_3-8\Gamma_4,\label{2.2}
\end{eqnarray}
where
\begin{equation}
\label{2.3} \ba{ll}
\ds\vs \Gamma_1:=b^{ij}b^m_ib_{jm}, & \ds
\Gamma_2:=b^{ij}b^m_k A^k_{\a j} A^\a_{mi}\\
\ds\vs \Gamma_3: =b^{ij}g_{ln}c^{\a\b}_{ji}c^{nl}_{\a\b},
& \Gamma_4:=b^j_lc^{\a\b}_{jk}c^{lk}_{\b\a}.\ea
\end{equation}
To continue, we need an expression for $\dt c^{\alpha\beta}_{ij}$.
\begin{eqnarray}
\ds \dt c^{\a\b}_{ij} & =&\ds\vs 2
a^{kl}A^\a_{ ik}A^\b_{lj}+ \dt A^\a_{ik} A^{\b k}_{j}
+ \dt A^\b_{kj} A^{\a k}_{i}\nonumber\\
&=&\ds\vs  2 a^{kl}A^\a _{ik}A^\b_{lj}\nonumber\\
&&\ds\vs + A^{\b k}_{ j}
(\D A^\a_{ik}-R^m_iA^\a _{mk}-R^m_k A^\a _{mi}+2A^\a_{mn}{{{R^m}_i}^n}_k)
-A^{\b k}_{ j}F^\a_m\l^m_{ik}\nonumber\\
&&\ds\vs + A^{\a k}_{ i}
(\D A^\b_{jk}-R^m_jA^\b _{mk}-R^m_k A^\b _{mj}+2A^\b_{mn}{{{R^m}_j}^n}_k)
-A^{\a k}_{ i}F^\b_m\l^m_{jk},\nonumber
\end{eqnarray}
where
\[\l^m_{ik}=\n_ib^m_k+\n_kb^m_i-\n^m b_{ik}.\]
Therefore
\begin{eqnarray}
\ds \dt c^{\a\b}_{ij}
&=&\ds\vs  2 a^{kl}A^\a _{ik}A^\b_{lj}+\D(A^\a_{ik}A^{\b k}_j)\nonumber\\
&&\ds\vs -\n ^l A^{\b k}_j\n_l A^\a_{ik}
+ A^{\b k}_{ j}(-R^m_iA^\a _{mk}-R^m_k A^\a _{mi}+2A^\a_{mn}{{{R^m}_i}^n}_k)
-A^{\b k}_{ j}F^\a_m\l^m_{ik}\nonumber\\
&&\ds\vs -\n ^l A^{\a k}_i\n_l A^\b_{jk} + A^{\a k}_{
i}(-R^m_jA^\b _{mk}-R^m_k A^\b _{mj}+2A^\b_{mn}{{{R^m}_j}^n}_k)
-A^{\a k}_{ i}F^\b_m\l^m_{jk}\nonumber\\
&=&
\Delta c^{\alpha\beta}_{ij}-2\n ^l A^{\b k}_j\n_l A^\a_{ik}
+2b^{kl}A^\a _{ik}A^\b_{lj}
-R_i^mc^{\alpha\beta}_{mj}-R_j^mc^{\alpha\beta}_{im}\nonumber\\
&&\ds\vs +2A^\alpha_{mn}A^{\beta k}_{j}\udud Rmink
+2A^{\alpha k}_{i}A^{\beta}_{mn}\udud Rmjnk
-A^{\b k}_{ j}F^\a_m\l^m_{ik}-A^{\a k}_{ i}F^\b_m\l^m_{jk}\label{3.1}
\end{eqnarray}

Then
\begin{eqnarray}
\dt\left(c^{\alpha\beta}_{ij}c^{ij}_{\beta\alpha}\right)
&=&2\dt\, c^{\alpha\beta}_{ij}c^{ij}_{\beta\alpha}
+4a_l^kc^{\alpha\beta}_{kj}c^{lj}_{\beta\alpha}\nonumber\\
&=&\Delta\left(c^{\alpha\beta}_{ij}c^{ij}_{\beta\alpha}\right)
-2\nabla^lc^{\alpha\beta}_{ij}\nabla_lc^{ij}_{\beta\alpha}
-4c^{ij}_{\beta\alpha}\nabla^lA^\alpha_{ik}\nabla_lA^{\beta k}_j\nonumber\\
&&+2c^{ij}_{\beta\alpha}\left\{2b^{kl}A^\a _{ik}A^\b_{lj}
-R_i^mc^{\alpha\beta}_{mj}-R_j^mc^{\alpha\beta}_{im}
+2A^\alpha_{mn}A^{\beta k}_{j}\udud Rmink
\right.\nonumber\\
&&\left.+2A^{\alpha k}_{i}A^{\beta}_{mn}\udud Rmjnk\right\}
+4a_l^kc^{\alpha\beta}_{kj}c^{lj}_{\beta\alpha}\nonumber\\
&=&\Delta\left(c^{\alpha\beta}_{ij}c^{ij}_{\beta\alpha}\right)
-2\nabla^lc^{\alpha\beta}_{ij}\nabla_lc^{ij}_{\beta\alpha}
-4c^{ij}_{\beta\alpha}\nabla^lA^\alpha_{ik}\nabla_lA^{\beta k}_j
\nonumber\\
&&+8b_l^kc^{\alpha\beta}_{kj}c^{lj}_{\beta\alpha}
+8c^{ij}_{\beta\alpha}A^\alpha_{mn}A^{\beta k}_j\udud Rmink\nonumber\\
&=&\Delta\left(c^{\alpha\beta}_{ij}c^{ij}_{\beta\alpha}\right)
-2\nabla^lc^{\alpha\beta}_{ij}\nabla_lc^{ij}_{\beta\alpha}
-4c^{ij}_{\beta\alpha}\nabla^lA^\alpha_{ik}\nabla_lA^{\beta k}_j
\nonumber\\
&&+8\Gamma_4+8\Gamma_5-8\Gamma_6,\label{3.2}
\end{eqnarray}
where
\begin{equation}\label{3.4}\ba{ll}
\ds\vs \Gamma_5:= c^{\alpha m}_{m\gamma}
c^{\beta\gamma}_{jl}c^{lj}_{\beta\alpha},
&  \ds
\Gamma_6:=c^{\a\gamma}_{mi}c^{\b m}_{j\gamma}c^{ij}_{\b\a}.
\ea\end{equation}
One can also show that
\begin{equation}\label{3.6}
\uudd R\alpha\beta ijR_{\a\b kl}R^{ijkl}=4(\Gamma_6-\Gamma_{7}),
\end{equation}
where
\[\Gamma_{7}:=
c^{\alpha j}_{i\beta}c^{l\gamma}_{\alpha j}c^{i\beta}_{\gamma l}.\]
Similarly, one has
\begin{equation}\label{3.7}
c^{\a m}_{m\gamma} {R_{\a\b}}^{ij}{R^{\gamma\b}}_{ij}=2(\Gamma_3-\Gamma_5),
\end{equation}
\begin{equation}\label{3.8}
\uudd R\alpha\beta ij\duud R\alpha\gamma il\dduu R\beta\gamma jl
=\Gamma_1-3\Gamma_4+3\Gamma_6-\Gamma_{7}
\end{equation}
and
\begin{equation}\label{3.9}
b^{ij}{R^{\a\b}}_{ki}{{R_{\a\b}}^k}_j=2(\Gamma_2-\Gamma_4).
\end{equation}
Let us define
\[G_1:=b^{mk}\n^lA^{\a}_{ik}\n_lA^i_{\a m}\,,\quad
G_2:=A^{k}_{\b j} A^i_{\a m} \n^l A^\a_{ik} \n_l A_{mj}^\b\]
and
\[G_3:=c^{mk}_{\a\b}\n^l A^\a_{ik}\n_l A^{\b i}_m\,,\quad
G_4:= A^{\b k}_jA^i_{\b m}\n^l A^\a_{ik}\n_l A_\a^{jm}.\]
Since
\[\n_l \uudd R\alpha\beta ij=\n_l A^\a_{ik} A^{\b k}_j+ A^\a_{ik} \n_l
A^{\b k}_j- \n_l A^\b_{ik} A^{\a k}_j- A^\b_{ik} \n_l A^{\a k}_j,\]
we obtain
\[\ba{rcl}
\ds\vs |\n \uudd R\alpha\beta ij|^2 &=&\ds
4 (b^{mk}\n A^\a_{ik}\n A^{i}_{\a m}+ A^k_{\beta j} A^i_{\a m} \n A^\a_{ik}
\n A_{mj}^\b\\
&& \ds\vs - c^{mk}_{\a\b}\n A^\a_{ik}A^{\b i}_m-A^{\b k}_jA^i_{\b m}\n
A^\a_{ik}\n A_\a^{jm})\\
&=&4(G_1+G_2-G_3-G_4).
\ea\]
In addition
\[\ba{rcl}
\ds\vs |\n_l b_{ij}|^2=2c^{kl}_{\a\b}\n A^\a_{ik}\n A^{\b i }_l+2
A^k_{\a j} A^{\b i}_l\n A^\a_{ik}\n A^{lj}_\b=2(G_2+G_3)
\ea\]
and
\[\ba{rcl}
\ds\vs \n(c^{\a\b}_{ij})\n(c_{\b\a}^{ij})
&=&2 A^\a_{ik} A^{jl}_\a\n A^{\b k}_j\n A^i_{\b l}
+2c^{\b l}_{k\a}\n A^{\a k}_i\n A^i_{\b l}\\
&=&2(G_3+G_4).
\ea\]
Combining (\ref{2.2}), (\ref{3.2}) and the last equations we get
\begin{eqnarray}
\dt|\uudd R\alpha\beta ij|^2&=&\Delta|\uudd R\alpha\beta ij|^2
-8(G_2+G_3)-8G_1+8(G_3+G_4)+8G_3\nonumber\\
&&+8(\Gamma_1+\Gamma_2+2\Gamma_3-2\Gamma_4)
-16(\Gamma_4+\Gamma_5-\Gamma_6)\nonumber\\
&=&\Delta|\uudd R\alpha\beta ij|^2-2|\nabla_l\uudd R\alpha\beta ij|^2
+8(\Gamma_1-3\Gamma_4+3\Gamma_6-\Gamma_7)
\nonumber\\
&&-8(\Gamma_6-\Gamma_7)
+16(\Gamma_3-\Gamma_5)+8(\Gamma_2-\Gamma_4)\nonumber
\end{eqnarray}
and in view of the previous equations for $\Gamma_i, i=1,\dots,7$,
we conclude
\begin{eqnarray}
\dt|\uudd R\alpha\beta ij|^2
&=&\Delta|\uudd R\alpha\beta ij|^2-2|\nabla_l\uudd R\alpha\beta ij|^2
+8\uudd R\alpha\beta ij\duud R\alpha\gamma il\dduu R\beta\gamma jl
\nonumber\\
&&-2R^{ijkl}\uudd R\alpha\beta ijR_{\a\b kl}
+8c^{\a m}_{m\gamma} {R_{\a\b}}^{ij}{R^{\gamma\b}}_{ij}
+4b^{ij}{R^{\a\b}}_{ki}{{R_{\a\b}}^k}_j\label{final}
\end{eqnarray}
and this equation is of the form given in (\ref{normal 2}).
\end{proof}
\begin{rem}
One can replace the gradient term $\nabla_l\uudd R\alpha\beta ij$ by
$\nabla_l^\perp\uudd R\alpha\beta ij,$ where $\nabla^\perp$ is the normal
connection. This simplifies (\ref{final}) a bit. We leave the details
to the reader.
\end{rem}

%%%%%%%%%%%%%%%%%%%%%%%%%%%%%%%%%%%%%%%%%%%%%%%%%%%%%%%%%%%%%%%%%%%%%%%%%%%%%%
\section{Graphical mean curvature flow with flat normal bundles}\label{graphs}
In the remaining
section we will see that in case of entire graphs with flat normal
bundles the computations in the paper by Ecker \& Huisken basically carry
over unchanged to the case of arbitrary codimension. In the following
we will outline the basic steps.

Let us assume that the initial submanifold
$M_0$ can be written as an entire graph over a flat $m$-plane in $\real{n}$
such that the normal bundle of $M_0$ is flat.

Suppose $\omega\in\Omega^m(\real{n})$ is a parallel $m$-form with
$|\omega|=1$. $\omega$ induces a function $w$ on any immersion
$F:M^m\to\real{n}$ by
$$F^*\omega=:wd\mu,$$
where $d\mu$ is the induced volume form on $M$.
Since $|\omega|=1$ we must have
$$-1\le w(x)\le 1\,,\quad\forall\,x\in M.$$
The angle $\alpha$ defined by $\cos\alpha=w$ measures the angle between
the flat plane defined by $\omega$ and the tangent planes of the submanifold
$M$. The condition to be a graph then is easily expressed by
$$w>0.$$
For hypersurfaces, $w$ is also given by the angle between the normal vector
$\nu$ of $M$ and a fixed normal direction of the reference plane defined by
$\omega$. This was considered in
\cite{ecker huisken}. For immersions with higher codimension,
a similar $w$-function was considered in \cite{jost xin}.
The evolution equation for $F^*\omega$ was derived earlier in a paper by
M.-T. Wang \cite{wang} and for the sake of completeness we include
the evolution equation in our notation.

\begin{equation}\label{omega}
\dt \omega_{i_1\cdots i_m}=\Delta\omega_{i_1\cdots i_m}
-R\omega_{i_1\cdots i_m}
-\sum_{s<j}\omega_{i_1\cdots i_{s-1}\alpha_si_{s+1}\cdots i_{j-1}
\alpha_ji_{j+1}\cdots i_m}\uudd R{\alpha_s}{\alpha_j}{i_s}{i_j},
\end{equation}
where $R$ denotes the scalar curvature and $\omega_{i_1\cdots i_m}
=\omega_{\a_1\cdots \a_m}F^{\a_1}_{i_1}\cdots F^{\a_m}_{i_m}$.

A first observation is
\begin{lem}\label{lem1}
Suppose $F:M^m\times[0,T)\to\real{n}$
is a smooth solution of the mean curvature flow
as in Theorem \ref{flatness}. Then $w$ defined
as above satisfies
\begin{equation}\label{w}
\dt w=\Delta w+w|A|^2.
\end{equation}
\end{lem}
\begin{proof}
We know that the normal bundle will be flat for all $t$. Hence
$$\dt F^*\omega=\Delta F^*\omega-RF^*\omega$$
and the conclusion follows from $\dt d\mu=-|H|^2d\mu$ and $R=|H|^2-|A|^2$.
\end{proof}

The second observation is
\begin{lem}\label{lem2}
Suppose $F:M^m\times[0,T)\to\real{n}$ is as in Lemma \ref{lem1}.
Then
\begin{eqnarray}
\dt |A|^2&=&\D |A|^2-2|\n^\perp A|^2+2|A^{\a m}_nA^{\b n}_m|^2\label{A2.1}\\
&\le&\D |A|^2-2|\n^\perp A|^2+2|A|^4\nonumber
\end{eqnarray}
\end{lem}
\begin{proof} The Lemma follows from (\ref{A2}), Theorem \ref{flatness}
and the estimate
$$|A^{\a m}_nA^{\b n}_m|^2\le|A|^4.$$
\end{proof}

Hence, in the case of a flat normal bundle the differential inequality
for $|A|^2$ is the same as in the case of codimension one.
Using these observations, we can follow \cite{ecker huisken} closely
to consider the mean curvature flow with flat normal bundles.

Now let us assume the following {\it linear growth} condition
\begin{equation}\label{condition}
v:=\frac{1}{w}\le c_1<\infty
\end{equation}
holds everywhere on $M_0$ for some constant $c_1$.
(\ref{w}) implies that if $M_0$ satisfies (\ref{condition})
then $M_t$ also satisfies (\ref{condition}) with the same constant.

\begin{lem}\label{lem3}The term $|A|^2v^2$ satisfies
\begin{equation}\label{Av}
\left(\dt -\D\right)|A|^2v^2\le -2v^{-1}\langle\n v,\n(|A|^2 v^2)\rangle.
\end{equation}
\end{lem}
\begin{proof}
To prove the Lemma, we need the following Kato type inequality
\begin{equation}\label{K}
|\n^\perp A|^2\ge \frac {n+2}n|\n|A||^2,
\end{equation}
which was proved in \cite{xin} for immersions with flat normal bundle.
Now the proof follows from the proof of Lemma 4.1 in \cite{ecker huisken}.
\end{proof}

Note that when the codimension is one, (\ref{K}) was proved in \cite{SSY}.
A direct consequence of Lemma \ref{lem3} and Corollary \ref{max princ} is
\begin{cor}\label{coro1}
If $M_t$ is a smooth solution of (\ref{mcf}) with bounded $v$
and bounded $|A|^2$ on each $M_t$, then there is the following estimate
\[ \sup_{M_t}|A|^2v^2 \le \sup_{M_0} |A|^2v^2.\]
\end{cor}
One can also get a differential inequality for higher derivatives
of the second fundmental form
\begin{eqnarray}
\ds \left(\dt -\D\right)\bigl(t^{l+1}|(\n^\perp)^l A|^2\bigr)
&\le&\ds\vs
 -2t^{l+1}|(\n^\perp)^{l+1} A|^2+(l+1)t^l|(\n^\perp)^l A|^2\nonumber\\
&&\hspace{-30pt}
+\, C(n,l)t^{l+1}\sum_{i+j+k=l}
|(\n^\perp)^i A||(\n^\perp)^j A||(\n^\perp)^l A||(\n^\perp)^l A|,\label{di}
\end{eqnarray}
for any integer $l\ge 0$, where $C(n,l)$ is a constant depending
only on $n$ and $l$. Therefore, the higher order estimates in
\cite{ecker huisken} can also be obtained in our case.
Furthermore we get
\begin{pro}\label{prop1}
Let $M_t$ be a smooth solution of (\ref{mcf}) with flat normal bundle.
Then for each $m\ge 0$ there is a consant $C(m)$ such
that
\begin{equation}\label{a1}
t^{m+1}|(\n^\perp)^mA|^2\le C(m).
\end{equation}
\end{pro}

Rescaling as in \cite{ecker huisken} (see also \cite{huisken2}),
we define
\[\widetilde F(s)=\frac 1{\sqrt{2t+1}} F(t),\]
where $s$ is given by $s=\frac 12 \log(2t+1)$. In the new time variable
$s$, we have a normalized equation
\begin{equation}\label{nmcf}
\frac{d}{ds} \widetilde F=\widetilde H-\widetilde F.
\end{equation}
{From} Proposition \ref{pro0}, Corollary \ref{coro1} and
Proposition \ref{prop1} we have estimates for the rescaled immersion
$\widetilde F$
\[\ba{rcl}
\ds \tilde u^2(\widetilde F,s)&\le&\ds\vs
\tilde c_0(1+\widetilde x^2
(\widetilde F,s))\\
\ds \widetilde v(\widetilde F,s)&\le &\ds\vs   c_1\\
\ds |\widetilde A|^2(\widetilde F,s) &\le& c_2,\ea\]
where $\tilde u^2:=|\widetilde F^\perp|^2$, $\tilde x^2=|\widetilde F^T|^2$
and $\tilde c_0$, $c_1$, $c_2$ are some
constants depending only on the initial immersion $F_0$.

Eventually, since the term $|\dd aij|^2$ appearing in the evolution
equation of the squared mean curvature (\ref{mean}) can be estimated by
$$|\dd aij|^2\le |H|^2|A|^2,$$
the computations can be carried out in the same way to derive
\begin{thm}\label{thm2} Suppose $M_0$
is an entire graph with bounded curvature over some $\real{m}$ in
$\real{n}$  satisfying the linear growth condition
(\ref{condition}) and that $M_0$ has a flat normal bundle.
Then the mean curvature flow admits a smooth solution for all $t>0$
with uniformly bounded curvature quantities.
If in addition we assume that
\begin{equation}\label{condition2}
u^2 \le c_3 (1+|F|^2)^{1-\d}
\end{equation}
holds on $M_0$ for some constant $c_3<\infty$ and $\d>0$, then the solution
$\widetilde M_s$ of the normalized mean curvature flow (\ref{nmcf})
converges for $s\to \infty$ to a
limiting surface $\widetilde M_\infty$ satisfying the equation for
selfsimilarly expanding solutions
\begin{equation}\label{solition}
F^\perp=H.
\end{equation}
\end{thm}

\begin{rem}
We note that Proposition 4.5 in \cite{ecker huisken},
the spatial decay behaviour, can also easily be done in the same way.
We leave the details to the reader.
\end{rem}

\begin{rem}
We improved the previous results of the third author. The dimension
limitation in \cite{xin} can be removed. That will appear in another paper.
\end{rem}
%%%%%%%%%%%%%%%%%%%%%%%%%%%%%%%%%%%%%%%%%%%%%%%%%%%%%%%%%%%%%%%%%%%%%%%%%%%%%%
\bibliographystyle{amsplain}

\end{document}